\documentclass[preprint]{elsarticle}

\usepackage{amssymb}
\usepackage{amsmath}
\usepackage{amsfonts}

\usepackage{hyperref,amstext}
\usepackage{mathrsfs}
\usepackage{longtable}
\usepackage{booktabs}

\usepackage[title]{appendix}
\usepackage{tikz}
\usepackage{subcaption}
\usepackage{algpseudocode, algorithm}

\newtheorem{theorem}{Theorem}

\newtheorem{lemma}{Lemma}

\newenvironment{proof}[1][{\rm PROOF.}]{\begin{trivlist}
\item[\hskip \labelsep {\bfseries #1}]}{\qed\end{trivlist}}

\journal{Discrete Applied Mathematics}

\begin{document}

\begin{frontmatter}



\title{An Improved Approximation for {\sc Maximum $k$-Dependent Set} on Bipartite Graphs}


\author{Seyedmohammadhossein Hosseinian, Sergiy Butenko}

\address{Texas A\&M University, Department of Industrial and Systems Engineering\\ College Station, TX 77843-3131, \{hosseinian, butenko\}@tamu.edu}

\begin{abstract}
We present a $(1+\frac{k}{k+2})$-approximation algorithm for the {\sc Maximum $k$-dependent Set} problem on bipartite graphs for any $k\ge1$. For a graph with $n$ vertices and $m$ edges, the algorithm runs in $\mathcal{O}(k m \sqrt{n})$ time and improves upon the previously best-known approximation ratio of $1+\frac{k}{k+1}$ established by Kumar et al.~[Theoretical Computer Science, 526: 90--96 (2014)]. Our proof also indicates that the algorithm retains its approximation ratio when applied to the (more general) class of K\"{o}nig-Egerv\'{a}ry graphs.
\end{abstract}

\begin{keyword}
maximum $k$-dependent set \sep approximation algorithms \sep bipartite graphs \sep {K\"{o}nig-Egerv\'{a}ry graphs} 


\end{keyword}

\end{frontmatter}




\section{Introduction} \label{sec:intro}
\noindent
Given a simple, undirected graph $G=(V,E)$, where $V$ is the set of vertices and $E$ is the set of edges, and a constant $k$, a subset of vertices $S\subseteq V$ is called a \emph{$k$-dependent set} if the degree of every vertex in the {subgraph induced by $S$} is at most $k$. The {\sc Maximum $k$-dependent Set} problem is to find a maximum-cardinality $k$-dependent set in $G$. This problem, first introduced in~\cite{DjidjevA92}, is a well-known generalization of the classical {\sc Maximum Independent Set} problem (the case of $k=0$) and is equivalent to the {\sc Maximum $s$-plex} problem~\cite{BBSBIVH11kplex,McCloskyHicks09} on the complement graph $\bar{G}$ of $G$, where $s=k+1$. For any $k \geq 0$, the {\sc Maximum $k$-dependent Set} problem is NP-hard on general graphs~\cite{DjidjevA92}, 
and by the theorem of Feige and Kogan~\cite{feige2005hardness}, it cannot be approximated (in polynomial time) within a factor of $n^{1-\epsilon}$, for any $\epsilon > 0$, unless NP=ZPP. Due to these hardness results, the {\sc Maximum $k$-dependent Set} problem has been also studied on special classes of graphs. It is known that, for any $k \geq 1$, the problem is polynomial-time solvable on cographs, split graphs, and graphs with bounded treewidth, hence trees, but it remains NP-hard on unit-disk graphs, planar graphs, and bipartite graphs~\cite{dessmarkA1993,havetA2009}. A constant-factor approximation algorithm and a polynomial time approximation scheme (PTAS) for this problem on unit-disk graphs have been presented in~\cite{balasundaramA2010} and~\cite{havetA2009}, respectively. Due to the Lipton-Tarjan separation theorem~\cite{liptonTarjan1979}, the problem is also known to admit PTAS on planar graphs. In this paper, we focus on the class of bipartite graphs. Given the NP-hardness result of Dessmark et al.~\cite{dessmarkA1993}, we are interested in solving the following problem approximately:
\begin{center}
\setlength{\fboxsep}{5pt}
\fbox{%
\parbox{114mm}{
{\sc Maximum $k$-Dependent Set in Bipartite Graphs (Max--$k$--DSBG)}
\medskip

\begin{tabular}{@{}ll}
{\bf Given:} & A bipartite graph $G=(V, E)$ and an integer $k\ge1$.\\
{\bf Find:} & A maximum-cardinality $k$-dependent set in $G$.
\end{tabular}
}%
}
\end{center}
We present a $(1+\frac{k}{k+2})$-approximation algorithm for {\sc Max--$k$--DSBG}, which is an improvement upon the best previously known approximation ratio of $1+\frac{k}{k+1}$, established by Kumar et al.~\cite{kumarA2014}.

To the best of our knowledge, the work of Kumar et al.~\cite{kumarA2014} is the only nontrivial approximation algorithm for {\sc Max--$k$--DSBG} in the literature. 
More specifically, Kumar et al.~\cite{kumarA2014} have studied approximability of the {\sc Maximum $\Pi_{k_1,k_2}$ Subgraph} and {\sc Minimum $\Pi_{k_1,k_2}$ Vertex Deletion} problems. Given a bipartite graph $G=(V_1\cup V_2, E)$ with the parts $V_1$ and $V_2$, let $\Pi_{k_1,k_2}$ be a property characterized by a finite set $\mathcal{H}$ of forbidden induced subgraphs, where any $H\in\mathcal{H}$ contains at most $k_1$ vertices in one of the parts $V_1$ and $V_2$ and at most $k_2$ vertices in the other of the parts. Then the {\sc Maximum $\Pi_{k_1,k_2}$ Subgraph} problem is to find a maximum-cardinality subset of vertices satisfying the property $\Pi_{k_1,k_2}$ (that is, containing no induced subgraph isomorphic to a graph from $\mathcal{H}$). The {\sc Minimum $\Pi_{k_1,k_2}$ Vertex Deletion} is the corresponding node deletion problem. Kumar et al.~\cite{kumarA2014} have shown that both these problems are APX-complete and proposed approximation algorithms based on the technique of iterative rounding applied to a sequence of optimal solutions of the linear programming relaxations. The approximation factors they established are given by $\max\{k_1,k_2\}$ for the {\sc Minimum $\Pi_{k_1,k_2}$ Vertex Deletion} and $2-\frac{1}{\max\{k_1,k_2\}}$ for the {\sc Maximum $\Pi_{k_1,k_2}$ Subgraph}, respectively.  It is evident that {\sc Max--$k$--DSBG} is equivalent to {\sc Maximum $\Pi_{k_1,k_2}$ Subgraph} if the corresponding {set $\mathcal H$ of forbidden subgraphs consists of a {\em star} graph $K_{1,k+1}$, i.e., a single vertex in one part adjacent to $k+1$ vertices in the other part of the graph}. Due to this equivalence, the method of Kumar et al.~\cite{kumarA2014} provides a $(1+\frac{k}{k+1})$-approximation algorithm for {\sc Max--$k$--DSBG}. Here, we present is a purely combinatorial approximation algorithm for this problem that achieves the approximation ratio of $1+\frac{k}{k+2}$ and runs in time bounded by $\mathcal{O}(k m \sqrt{n})$ on a graph with $n$ vertices and $m$ edges. Our method is surprisingly simple. 

We use the following concepts and notations.  Let $n := |V|$ and $m:=|E|$ denote the number of vertices and edges in $G$, respectively. For a subset of vertices $S$, $G[S]$ denotes the subgraph induced by $S$ in $G$. A subset of edges $M\subseteq E$ is called a matching if no two edges in $M$ have a common incident vertex. A subset of vertices $C\subseteq V$ is a vertex cover if every edge in $E$ has at least one of its end-vertices in $C$. K\"{o}nig-Egerv\'{a}ry matching theorem~\cite{denes1931grafok}  (see
\cite{Rizzi} for a short proof) states that 
the maximum cardinality of a matching in 
a bipartite graph $G$ is equal to the minimum cardinality of a vertex cover.  The reader is referred to~\cite{Ahuja} and~\cite{Vazirani03} for background on polynomial-time  algorithms for problems in graphs/networks and approximation algorithms, respectively.

\section{The Algorithm} \label{sec:alg}
\noindent
In this section, we present our algorithm, establish its approximation ratio, and analyze its time complexity. 

The proposed method consists of two simple stages. The first stage is comprised of $k$ recursive iterations, each of which computes a maximum-cardinality matching and removes it from the graph. This is followed by the second stage, which computes a maximum independent set $S$ in the residual graph and outputs it as an approximate solution for {\sc Max--$k$--DSBG}. The process is formally described in Algorithm~\ref{algorithm}. 
Notice~that, in line 3 of the algorithm, some {\em edges} of the graph are deleted, not their {\em end-vertices}. An illustration of the process for $k=1$ on an example graph is provided in Figure~\ref{fig:k_one}. 

\begin{algorithm}[t]
\caption{\sc An approximation algorithm for {\sc Max--$k$--DSBG}}\label{algorithm} 
\begin{algorithmic}[1]
\Require{An integer constant $k\ge1$ and a bipartite graph $G=(V, E)$}
\medskip

\For{$i=1$ to $k$}
\State \quad find a maximum-cardinality matching $M \subseteq E$ in $G=(V,E)$
\State \quad $E\gets E\setminus M$ 
\EndFor
\State find a maximum independent set $S$ in the residual graph $G$
\State\Return{$S$}
\end{algorithmic}
\end{algorithm}

\tikzstyle{label}=[fill=black!0,minimum size=15pt,inner sep=0pt]
\tikzstyle{vertex}=[draw, circle,fill=black!0,minimum size=15pt,inner sep=0pt]
\tikzstyle{darkVertex}=[draw,circle,fill=black!25,minimum size=15pt,inner sep=0pt]
\tikzstyle{voidVertex}=[draw,white,circle,fill=black!0,minimum size=15pt,inner sep=0pt]
\tikzstyle{edge} = [draw,thick,-]
\tikzstyle{thinEdge} = [draw,thin,-]
\tikzstyle{semiEdge} = [draw,semithick,-]
\tikzstyle{ultraEdge} = [draw,ultra thick,-]
\tikzstyle{dashEdge} = [draw,dashed,-]
\tikzstyle{dotEdge} = [draw,dotted,-]
\tikzstyle{weight} = [font=\small]
\begin{figure}[ht]
\begin{subfigure}{.5\textwidth}
\centering
\begin{tikzpicture}[scale=0.5, auto,swap]
\foreach \pos/\name in {{(0,15)/1}, {(4,15)/2}, {(0,12)/3}, {(4,12)/4}, {(4,9)/6}, {(0,6)/7}, {(4,6)/8}, {(0,3)/9}, {(4,3)/10}}
\node[vertex] (\name) at \pos {$\name$};
\foreach \pos/\name in {{(0,9)/5}, {(0,0)/11}, {(4,0)/12}}
\node[darkVertex] (\name) at \pos {$\name$};
\foreach \source/ \dest /\weight in {1/2/, 1/12/, 2/5/, 2/11/, 3/4/, 3/12/, 5/6/, 6/11/, 7/8/, 7/12/, 9/10/, 10/11/}
\path[edge] (\source) -- node[weight]{}(\dest);
\end{tikzpicture}
\caption{}
\label{fig:1a}
\end{subfigure}%
\begin{subfigure}{.5\textwidth}
\centering
\begin{tikzpicture}[scale=0.5, auto,swap]
\foreach \pos/\name in {{(0,15)/1}, {(0,12)/3}, {(4,12)/4}, {(4,9)/6}, {(0,6)/7}, {(4,6)/8}, {(0,3)/9}, {(4,3)/10}}
\node[vertex] (\name) at \pos {$\name$};
\foreach \pos/\name in {{(4,15)/2}, {(0,9)/5}, {(0,0)/11}, {(4,0)/12}}
\node[darkVertex] (\name) at \pos {$\name$};
\foreach \source/ \dest /\weight in {1/2/, 2/11/, 3/12/, 5/6/, 7/12/, 10/11/}
\path[edge] (\source) -- node[weight]{}(\dest);
\foreach \source/ \dest /\weight in {1/12/, 2/5/, 3/4/, 6/11/, 7/8/, 9/10/}
\path[dashEdge] (\source) -- node[weight]{}(\dest);
\end{tikzpicture}
\caption{}
\label{fig:1b}
\end{subfigure}
\caption{An example for $k=1$: (a) the set of white vertices is a maximum 1-dependent set, (b) a maximum-cardinality matching $M$ is illustrated by the dashed edges, and the gray vertices depict a minimum vertex cover of the residual graph after deleting $M$; the set of white vertices is the output of the algorithm.}
\label{fig:k_one}
\end{figure}
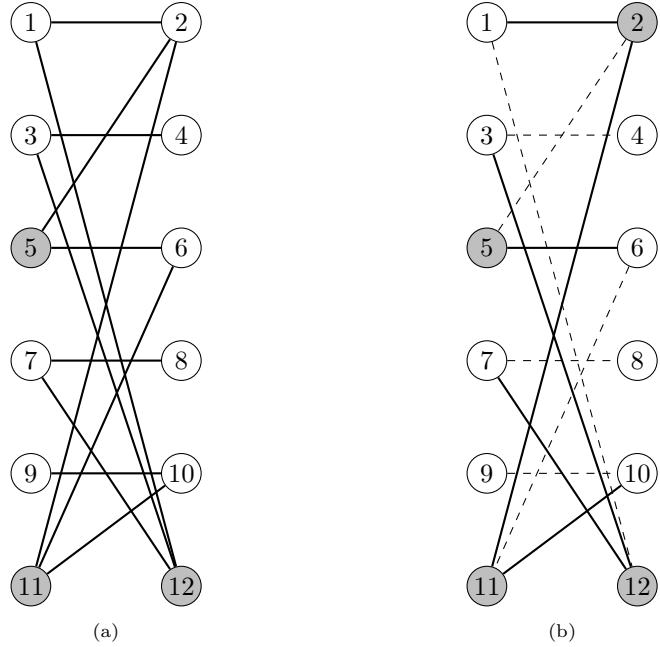

To establish the correctness of the algorithm, observe that a vertex $v \in S$ can be incident to at most $k$ (deleted) edges in $G[S]$; hence, $S$ is a feasible solution to {\sc Max--$k$--DSBG}. We proceed to establish the approximation ratio of the algorithm. 

Let $\mathcal{M}_i,~i\in \{1,\ldots,k\}$, be the {\em accumulative} set of edges deleted through the first $i$ iterations of lines 2 and 3 of Algorithm~\ref{algorithm}, and let $G_i$ denote the corresponding residual graph, i.e., the graph obtained from $G$ by removing $\mathcal{M}_i$. Furthermore, let $I_i,~i\in \{1,\ldots,k\}$, be a maximum independent set of $G_i$, and let $I^*_k$ be an arbitrary maximum $k$-dependent set in $G$; by this notation, $S \equiv I_k$. 
Patently, the edge set of the subgraph induced by $I^*_k$, i.e., $E(G[I^*_k])$, can be decomposed into a collection of disjoint matchings in $G$; removing the set of edges belonging to these matchings from the graph, and then, finding a maximum independent set in the residual graph generates an optimal solution to {\sc Max--$k$--DSBG} on $G$. The intuition behind our proof is based on observation that the performance of Algorithm~\ref{algorithm} depends on the number of edges from $E(G[I^*_k])$ that are not included in the maximum-cardinality matchings detected through $k$ iterations of the algorithm.

\begin{lemma} \label{lem:matching}
$|I_k| \geq n - \dfrac{|\mathcal{M}_k|}{k}$, for every fixed $k \geq 1$.
\end{lemma}

\begin{proof}
For $k=1$,
$$
|I_1| \geq \alpha (G) = n - |\mathcal{M}_1|,
$$
where $\alpha(G)$ denotes the independence number of $G$, and the equality is due to the K\"{o}nig-Egerv\'{a}ry matching theorem. For $k \geq 2$, the proof is by induction. The base case is already established, so for a fixed $k \geq 2$ and some $i \in \{2,\ldots,k\}$, suppose that
$$
|I_{i-1}| \geq n - \dfrac{|\mathcal{M}_{i-1}|}{i-1}.
$$
Then, using the induction hypothesis and the relation $|M(G_{i-1})| = n - |I_{i-1}|$ given by the K\"{o}nig-Egerv\'{a}ry matching theorem, we obtain
$$
|\mathcal{M}_i| = |\mathcal{M}_{i-1}|+|M(G_{i-1})| \geq (i-1) \big ( n - |I_{i-1}| \big ) + \big ( n - |I_{i-1}| \big ) = i \big ( n - |I_{i-1}| \big ),
$$
where $M(G_{i-1})$ denotes a maximum-cardinality matching in $G_{i-1}$. 
{Finally, the fact that $G_i$ is obtained by removing $M(G_{i-1})$ from $G_{i-1}$, implies that $|I_{i-1}|\le|I_i|$; hence, 
$$
|\mathcal{M}_i|\ge i \big ( n - |I_{i}| \big ).
$$}
This completes the proof.
\end{proof}

\bigskip

\noindent
Let $E'_k := {E(G[I_k^*])} \backslash \mathcal{M}_k$ denote the set of edges with both end-vertices in the optimal solution $I_k^*$ that are not deleted through executions of lines 2 and 3 of Algorithm~\ref{algorithm}. 

\begin{lemma} \label{lem:main}
$|I_k| \geq \dfrac{|I^*_k|}{2} + \dfrac{|E'_k|}{k}$, for every fixed $k \geq 1$.
\end{lemma}

\begin{proof}
Let $D^*_k := V \backslash I^*_k$ be the set of vertices not included in the  optimal solution $I^*_k$, and $E''_k := {E(G[I^*_k])} \backslash E'_k = {E(G[I^*_k])} \cap \mathcal{M}_k$. Then,
$$
|\mathcal{M}_k| \leq |E''_k| + k |D^*_k|  = \big ( |{E(G[I^*_k])}| - |E'_k| \big ) + k |D^*_k| \leq \dfrac{k|I^*_k|}{2} - |E'_k| + k |D^*_k|,
$$
where the first inequality is due to the fact that every edge from $\mathcal{M}_k$ not belonging to $E''_k$ has at least one end-vertex in $D^*_k$, and at most $k$ edges containing a vertex $v \in D^*_k$ can belong to $\mathcal{M}_k$, i.e., at most one edge from each of the $k$ computed matchings.  The second inequality is implied by the fact that the degrees of the vertices of $I^*_k$ in $G[I^*_k]$ are bounded by $k$. This result, along with Lemma~\ref{lem:matching}, leads to
$$
|I_k| \geq n - \dfrac{|\mathcal{M}_k|}{k} \geq \big ( n - |D^*_k| \big ) - \dfrac{|I^*_k|}{2} + \dfrac{|E'_k|}{k} = |I^*_k| - \dfrac{|I^*_k|}{2} + \dfrac{|E'_k|}{k} = \dfrac{|I^*_k|}{2} + \dfrac{|E'_k|}{k}.
$$
\end{proof}


\begin{lemma} \label{lem:ratio}
$\dfrac{|I_k|}{|I^*_k|} \geq \dfrac{1}{2} + \dfrac{1}{2(k+1)}$, for every fixed $k \geq 1$.
\end{lemma}

\begin{proof}
Suppose $|I^*_k| \leq \dfrac{2|E'_k|(k+1)}{k}$, or equivalently $\dfrac{|E'_k|}{k} \geq \dfrac{|I^*_k|}{2(k+1)}$. This, together with Lemma~\ref{lem:main}, implies
$$
|I_k| \geq \dfrac{|I^*_k|}{2} + \dfrac{|E'_k|}{k} \geq \dfrac{|I^*_k|}{2} + \dfrac{|I^*_k|}{2(k+1)} = |I^*_k| \left ( \dfrac{1}{2} + \dfrac{1}{2(k+1)} \right ).
$$
Now, consider the case $|I^*_k| > \dfrac{2|E'_k|(k+1)}{k}$, or equivalently $|E'_k| < \dfrac{k|I^*_k|}{2(k+1)}$. Observe that, at most $|E'_k|$ vertices (one for each edge in $G_k[I_k^*]$) should be deleted from $I^*_k$ to turn it into an independent set in $G_k$. Thus, $|I_k| \geq |I^*_k| - |E'_k|$, which, under the case assumption, leads to
$$
|I_k| \geq |I^*_k| - |E'_k| > |I^*_k| - \dfrac{k|I^*_k|}{2(k+1)} = |I^*_k| \left ( 1- \dfrac{k}{2(k+1)} \right ) = |I^*_k| \left ( \dfrac{1}{2} + \dfrac{1}{2(k+1)} \right ).
$$
\end{proof}

\begin{theorem}
Algorithm~\ref{algorithm} is a $(1+\frac{k}{k+2})$-approximation algorithm for {\sc Max--$k$--DSBG} that runs in $\mathcal{O}(km \sqrt{n})$ time.
\end{theorem}

\begin{proof}
By Lemma~\ref{lem:ratio} and the feasibility of an output of Algorithm~\ref{algorithm} to {\sc Max--$k$--DSBG},
$$
\dfrac{|I_k^*|}{|S|}=\dfrac{|I_k^*|}{|I_k|}\le 1+\frac{k}{k+2},
$$
which establishes the approximation ratio. 

Finding a maximum-cardinality matching $M$ in a bipartite graph $G=(V,E)$ takes no more than $\mathcal{O}(m \sqrt{n})$ time by the Hopcroft-Karp algorithm~\cite{hopcroftA1973}. Deleting $M$ itself takes no more than a linear time. Thus, the loop of Algorithm~\ref{algorithm} runs in a time bounded by $\mathcal{O}(km \sqrt{n})$. Besides, by the {K\"{o}nig-Egerv\'{a}ry} matching theorem, the Hopcroft-Karp algorithm can be also used to find a maximum independent set in the residual graph $G_k$. Therefore, the time complexity of Algorithm~\ref{algorithm} is no worse than $\mathcal{O}(km \sqrt{n})$.
\end{proof}

{
It is worth noting that 1-dependent set and {\em independent union of cliques} (IUC) are equivalent in bipartite graphs, so Algorithm~\ref{algorithm} (for $k=1$) provides a $\frac{4}{3}$-approximation algorithm for the {\sc Maximum IUC} problem~\cite{ertem2020IUC,hosseinian2021polyhedral} on this class of graphs. We should also point out that finding a maximum-cardinality matching can be done by randomized algorithms with a better time complexity than the Hopcroft-Karp method; see e.g.~\cite{mucha2004maximum}.
}

{
Bipartite graphs are a subclass of K\"{o}nig-Egerv\'{a}ry graphs. A graph is called a K\"{o}nig-Egerv\'{a}ry (KE) graph if the maximum cardinality of a matching is equal to the minimum cardinality of a vertex cover in the graph. It is known that the {\sc Maximum Independent Set} problem is polynomial-time solvable on KE graphs~\cite{deming1979,sterboul1979}. Given this result, as well as the fact that our proof relies only on the K\"{o}nig-Egerv\'{a}ry matching theorem, Algorithm~\ref{algorithm} is indeed a $(1+\frac{k}{k+2})$-approximation algorithm for the {\sc Maximum $k$-dependent Set} problem on the (more general) class of KE graphs. A forbidden-subgraph characterization of KE graphs is presented in~\cite{korachA2006}.
}

\section{Examples Achieving the Approximation Ratio Bound}
\label{sec:tight_example}
\noindent
In this section, we present a family of bipartite graphs, on which Algorithm~\ref{algorithm} may achieve the  $1+\frac{k}{k+2}=\frac{2(k+1)}{k+2}$ approximation bound.

For a fixed $k \geq 1$, let $G=(V_1 \cup V_2, E)$ be a bipartite graph on $2(k+2)$ vertices, where $V_1 = A \cup \{u\}$, $V_2 = B \cup \{w\}$, $A=\{a_1,\ldots, a_{k+1}\}$, and $B=\{b_1,\ldots, b_{k+1}\}$. We distinguish $u$ and $w$ from the rest of the vertices, as they are the only ones not included in the optimal solution of {\sc Max--$k$--DSBG} on $G$. The construction of $G$ is as follows: $u$ and $w$ are adjacent; every vertex in $A$ (resp. $B$) is adjacent to $w$ (resp. $u$) as well as all vertices in $B$ (resp. $A$) except for the one with the same index. That is,
\begin{equation*}
\begin{aligned}
E= &\{\{u,w\}\}\cup\{\{a_i,b_{i'}\}\mid i \in \{1,\ldots,k+1\},~i' \in \{1,\ldots,k+1\} \backslash \{i\}\}\\
&\cup\{\{a_i,w\}\mid i \in \{1,\ldots,k+1\}\}\cup \{\{b_i,u\}\mid i \in \{1,\ldots,k+1\}\}. 
\end{aligned}
\end{equation*}
Observe that, by this construction, $A \cup B$ is the (unique) optimal solution of {\sc Max--$k$--DSBG} on $G$, with the optimal solution value $|A \cup B|=2(k+1)$. Besides, $G[A \cup B]$ has $k$ disjoint perfect matchings, each of which constitutes a perfect matching for $G$ with $\{u,w\} \in E$.

We show that the output of Algorithm~\ref{algorithm} on $G$ may contain exactly $k+2$ vertices; we present the process for even and odd values of $k$ separately.

Let $k$ be even, and consider the performance of Algorithm~\ref{algorithm} on $G$ as follows: at every iteration $j \in \{1,\ldots,k\}$, the algorithm removes the (perfect) matching 
\begin{equation} \label{eq:match1}
M_j=(\{\{a_i , b_{i+j}\}\mid i \in \{1,..,k+1\}\} \backslash \{a_j,b_{2j}\}) \cup \{ \{a_j,w\},\{b_{2j},u\} \}
\end{equation}
from $G$, with the convention $b_{i+k+1} \equiv b_i$. Note that, at every iteration $j \in \{1,\ldots,k\}$, the algorithm includes the edges $\{a_j,w\}$ and $\{b_{2j},u\}$ in $M_j$ and leaves the edges $\{a_j,b_{2j}\}$ and $\{u,w\}$ intact. As a result, after $k$ iterations of the loop of Algorithm~\ref{algorithm}, the set of edges $\widetilde{E}=\{\{a_j,b_{2j}\},~\forall j \in \{1,\ldots,k\}\}$ will remain in the residual graph, i.e., the graph obtained from $G$ by removing $\mathcal{M}_k = \bigcup_{j=1}^k M_j$. Observe that, given $k$ is even, there exists no index $j \in \{1,\ldots,k\}$ such that $2j=2j'+k+1$ for some index $j' \in \{1,\ldots,j-1\}$. This implies that no two edges in $\widetilde{E}$ share an end-vertex, so every minimum vertex cover of the residual graph contains one end-vertex for each edge in $\widetilde{E}$. Besides, it is easy to verify that, in addition to $\widetilde{E}$, the edge set of the residual graph contains (exactly) three edges, i.e., $\{a_{k+1},w\}$, $\{w,u\}$, and $\{u,b_{k+1}\}$, which form the (isolated) path $P=\{a_{k+1},w,u,b_{k+1}\}$. Every minimum vertex cover of this graph must contain two out of four vertices from $P$. This implies that every independent set of the residual graph is of cardinality $2(k+2)-(k+2)=k+2$.

Now, let $k$ be odd, and consider the abovementioned process with a slight modification: at every iteration $j \in \{1,\ldots, \frac{k+1}{2}\}$, the algorithm removes the matching defined by~\eqref{eq:match1}, and at every iteration $j \in \{\frac{k+3}{2},\ldots, k\}$, it removes the following (perfect) matching from $G$:
$$
M_j=(\{\{a_i, b_{i+j}\},~\forall i \in \{1,..,k+1\}\} \backslash \{a_{j+1},b_{2j+1}\}) \cup \{ \{a_{j+1},w\},\{b_{2j+1},u\} \}.
$$
Notice that $b_{k+3} \equiv b_2$ for every odd value of $k$; at the iteration $j=\frac{k+3}{2}$, instead of the edge $\{a_{\frac{k+3}{2}},b_{2}\}$, the algorithm excludes the edge $\{a_{\frac{k+3}{2}+1},b_{3}\}$ from the matching $M_j$, and it follows this pattern for the rest of the iterations. This ensures that the set of edges connecting $A$ and $B$ in the residual graph, denoted by $\widetilde{E}$, are isolated, i.e., no two edges among them have a common end-vertex. Besides, similar to the previous case, the residual graph always includes the path $P=\{a_{\frac{k+3}{2}},w,u,b_1\}$, which indicates that every independent set of this graph is of cardinality $2(k+2)-(k+2)=k+2$.

\begin{figure}[!ht]
\begin{subfigure}{.33\textwidth}
\centering
\begin{tikzpicture}[scale=0.45,auto,swap]
\foreach \pos/\name in {{(1,16)/A},{(5,16)/B},{(0,15)/1},{(0,12)/2},{(0,9)/3},{(0,6)/4},{(0,3)/5}}
\node[label] (\name) at \pos {$\name$};
\foreach \pos/\name in {{(1,15)/A1}, {(5,15)/B1}, {(1,12)/A2}, {(5,12)/B2}, {(1,9)/A3}, {(5,9)/B3}, {(1,6)/A4}, {(5,6)/B4}, {(1,3)/A5}, {(5,3)/B5}}
\node[vertex] (\name) at \pos {};
\foreach \pos/\name in {{(1,0)/u}, {(5,0)/w}}
\node[darkVertex] (\name) at \pos {$\name$};
\foreach \source/ \dest /\weight in {A1/B2/, A1/B3/, A1/B4/, A1/B5/, A2/B1/, A2/B3/, A2/B4/, A2/B5/, A3/B1/, A3/B2/, A3/B4/, A3/B5/, A4/B1/, A4/B2/, A4/B3/, A4/B5/, A5/B1/, A5/B2/, A5/B3/, A5/B4/}
\path[semiEdge] (\source) -- node[weight]{}(\dest);
\foreach \source/ \dest /\weight in {A1/w/, A2/w/, A3/w/, A4/w/, A5/w/, B1/u/, B2/u/, B3/u/, B4/u/, B5/u/, u/w/}
\path[semiEdge] (\source) -- node[weight]{}(\dest);
\end{tikzpicture}
\caption{}
\label{fig:1ex_4}
\end{subfigure}
\begin{subfigure}{.33\textwidth}
\centering
\begin{tikzpicture}[scale=0.45,auto,swap]
\foreach \pos/\name in {{(1,16)/A},{(5,16)/B},{(0,15)/1},{(0,12)/2},{(0,9)/3},{(0,6)/4},{(0,3)/5}}
\node[label] (\name) at \pos {$\name$};
\foreach \pos/\name in {{(1,15)/A1}, {(5,15)/B1}, {(1,12)/A2}, {(5,12)/B2}, {(1,9)/A3}, {(5,9)/B3}, {(1,6)/A4}, {(5,6)/B4}, {(1,3)/A5}, {(5,3)/B5}}
\node[vertex] (\name) at \pos {};
\foreach \pos/\name in {{(1,0)/u}, {(5,0)/w}}
\node[vertex] (\name) at \pos {$\name$};
\foreach \source/ \dest /\weight in {A1/B3/, A1/B4/, A1/B5/, A2/B1/, A2/B4/, A2/B5/, A3/B1/, A3/B2/, A3/B5/, A4/B1/, A4/B2/, A4/B3/, A5/B2/, A5/B3/, A5/B4/}
\path[semiEdge] (\source) -- node[weight]{}(\dest);
\foreach \source/ \dest /\weight in {A2/w/, A3/w/, A4/w/, A5/w/, B1/u/, B3/u/, B4/u/, B5/u/, u/w/}
\path[semiEdge] (\source) -- node[weight]{}(\dest);
\foreach \source/ \dest /\weight in {A2/B3/, A3/B4/, A4/B5/, A5/B1/}
\path[dashEdge] (\source) -- node[weight]{}(\dest);
\foreach \source/ \dest /\weight in {A1/B2/}
\path[ultraEdge] (\source) -- node[weight]{}(\dest);
\foreach \source/ \dest /\weight in {A1/w/, B2/u/}
\path[dotEdge] (\source) -- node[weight]{}(\dest);
\end{tikzpicture}
\caption{}
\label{fig:2ex_4}
\end{subfigure}
\begin{subfigure}{.32\textwidth}
\centering
\begin{tikzpicture}[scale=0.45,auto,swap]
\foreach \pos/\name in {{(1,16)/A},{(5,16)/B},{(0,15)/1},{(0,12)/2},{(0,9)/3},{(0,6)/4},{(0,3)/5}}
\node[label] (\name) at \pos {$\name$};
\foreach \pos/\name in {{(1,15)/A1}, {(5,15)/B1}, {(1,12)/A2}, {(5,12)/B2}, {(1,9)/A3}, {(5,9)/B3}, {(1,6)/A4}, {(5,6)/B4}, {(1,3)/A5}, {(5,3)/B5}}
\node[vertex] (\name) at \pos {};
\foreach \pos/\name in {{(1,0)/u}, {(5,0)/w}}
\node[vertex] (\name) at \pos {$\name$};
\foreach \source/ \dest /\weight in {A1/B4/, A1/B5/, A2/B1/, A2/B5/, A3/B1/, A3/B2/, A4/B2/, A4/B3/, A5/B3/, A5/B4/}
\path[semiEdge] (\source) -- node[weight]{}(\dest);
\foreach \source/ \dest /\weight in {A3/w/, A4/w/, A5/w/, B1/u/, B3/u/, B5/u/, u/w/}
\path[semiEdge] (\source) -- node[weight]{}(\dest);
\foreach \source/ \dest /\weight in {A2/B3/, A3/B4/, A4/B5/, A5/B1/, A1/B3/, A3/B5/, A4/B1/, A5/B2/}
\path[dashEdge] (\source) -- node[weight]{}(\dest);
\foreach \source/ \dest /\weight in {A1/B2/, A2/B4/}
\path[ultraEdge] (\source) -- node[weight]{}(\dest);
\foreach \source/ \dest /\weight in {A1/w/, B2/u/, A2/w/, B4/u/}
\path[dotEdge] (\source) -- node[weight]{}(\dest);
\end{tikzpicture}
\caption{}
\label{fig:3ex_4}
\end{subfigure}
\begin{subfigure}{.33\textwidth}
\centering
\begin{tikzpicture}[scale=0.45,auto,swap]
\foreach \pos/\name in {{(1,16)/A},{(5,16)/B},{(0,15)/1},{(0,12)/2},{(0,9)/3},{(0,6)/4},{(0,3)/5}}
\node[label] (\name) at \pos {$\name$};
\foreach \pos/\name in {{(1,15)/A1}, {(5,15)/B1}, {(1,12)/A2}, {(5,12)/B2}, {(1,9)/A3}, {(5,9)/B3}, {(1,6)/A4}, {(5,6)/B4}, {(1,3)/A5}, {(5,3)/B5}}
\node[vertex] (\name) at \pos {};
\foreach \pos/\name in {{(1,0)/u}, {(5,0)/w}}
\node[vertex] (\name) at \pos {$\name$};
\foreach \source/ \dest /\weight in {A1/B5/, A2/B1/, A3/B2/, A4/B3/, A5/B4/}
\path[semiEdge] (\source) -- node[weight]{}(\dest);
\foreach \source/ \dest /\weight in {A4/w/, A5/w/, B3/u/, B5/u/, u/w/}
\path[semiEdge] (\source) -- node[weight]{}(\dest);
\foreach \source/ \dest /\weight in {A2/B3/, A3/B4/, A4/B5/, A5/B1/, A1/B3/, A3/B5/, A4/B1/, A5/B2/, A1/B4/, A2/B5/, A4/B2/, A5/B3/}
\path[dashEdge] (\source) -- node[weight]{}(\dest);
\foreach \source/ \dest /\weight in {A1/B2/, A2/B4/, A3/B1/}
\path[ultraEdge] (\source) -- node[weight]{}(\dest);
\foreach \source/ \dest /\weight in {A1/w/, B2/u/, A2/w/, B4/u/, A3/w/, B1/u/}
\path[dotEdge] (\source) -- node[weight]{}(\dest);
\end{tikzpicture}
\caption{}
\label{fig:4ex_4}
\end{subfigure}
\begin{subfigure}{.33\textwidth}
\centering
\begin{tikzpicture}[scale=0.45,auto,swap]
\foreach \pos/\name in {{(1,16)/A},{(5,16)/B},{(0,15)/1},{(0,12)/2},{(0,9)/3},{(0,6)/4},{(0,3)/5}}
\node[label] (\name) at \pos {$\name$};
\foreach \pos/\name in {{(1,15)/A1}, {(5,15)/B1}, {(1,12)/A2}, {(5,12)/B2}, {(1,9)/A3}, {(5,9)/B3}, {(1,6)/A4}, {(5,6)/B4}, {(1,3)/A5}, {(5,3)/B5}}
\node[vertex] (\name) at \pos {};
\foreach \pos/\name in {{(1,0)/u}, {(5,0)/w}}
\node[vertex] (\name) at \pos {$\name$};
\foreach \source/ \dest /\weight in {}
\path[semiEdge] (\source) -- node[weight]{}(\dest);
\foreach \source/ \dest /\weight in {A5/w/, B5/u/, u/w/}
\path[semiEdge] (\source) -- node[weight]{}(\dest);
\foreach \source/ \dest /\weight in {A2/B3/, A3/B4/, A4/B5/, A5/B1/, A1/B3/, A3/B5/, A4/B1/, A5/B2/, A1/B4/, A2/B5/, A4/B2/, A5/B3/, A1/B5/, A2/B1/, A3/B2/, A5/B4/}
\path[dashEdge] (\source) -- node[weight]{}(\dest);
\foreach \source/ \dest /\weight in {A1/B2/, A2/B4/, A3/B1/, A4/B3/}
\path[ultraEdge] (\source) -- node[weight]{}(\dest);
\foreach \source/ \dest /\weight in {A1/w/, B2/u/, A2/w/, B4/u/, A3/w/, B1/u/, A4/w/, B3/u/}
\path[dotEdge] (\source) -- node[weight]{}(\dest);
\end{tikzpicture}
\caption{}
\label{fig:5ex_4}
\end{subfigure}
\begin{subfigure}{.32\textwidth}
\centering
\begin{tikzpicture}[scale=0.45,auto,swap]
\foreach \pos/\name in {{(1,16)/A},{(5,16)/B},{(0,15)/1},{(0,12)/2},{(0,9)/3},{(0,6)/4},{(0,3)/5}}
\node[label] (\name) at \pos {$\name$};
\foreach \pos/\name in {{(1,15)/A1}, {(1,12)/A2}, {(1,9)/A3}, {(1,6)/A4}, {(1,3)/A5}, {(5,3)/B5}}
\node[vertex] (\name) at \pos {};
\foreach \pos/\name in {{(1,0)/u}, {(5,0)/w}}
\node[darkVertex] (\name) at \pos {$\name$};
\foreach \pos/\name in {{(5,15)/B1}, {(5,12)/B2}, {(5,9)/B3}, {(5,6)/B4}}
\node[darkVertex] (\name) at \pos {};
\foreach \source/ \dest /\weight in {A1/B2/, A1/B3/, A1/B4/, A1/B5/, A2/B1/, A2/B3/, A2/B4/, A2/B5/, A3/B1/, A3/B2/, A3/B4/, A3/B5/, A4/B1/, A4/B2/, A4/B3/, A4/B5/, A5/B1/, A5/B2/, A5/B3/, A5/B4/}
\path[semiEdge] (\source) -- node[weight]{}(\dest);
\foreach \source/ \dest /\weight in {A1/w/, A2/w/, A3/w/, A4/w/, A5/w/, B1/u/, B2/u/, B3/u/, B4/u/, B5/u/, u/w/}
\path[semiEdge] (\source) -- node[weight]{}(\dest);
\end{tikzpicture}
\caption{}
\label{fig:6ex_4}
\end{subfigure}
\caption{Worst-case example for $k=4$: (a) the original graph; the set of white vertices is a maximum 4-dependent set, (b) the graph after removing $\mathcal{M}_1$, (c) the graph after removing $\mathcal{M}_2$, (d) the graph after removing $\mathcal{M}_3$, (e) the graph after removing $\mathcal{M}_4$, (f) the original graph; the gray vertices depict a minimum vertex cover of the residual graph, i.e., after removing $\mathcal{M}_4$, and the set of white vertices is the output of the algorithm.}
\label{fig:tight_ex_4}
\end{figure}
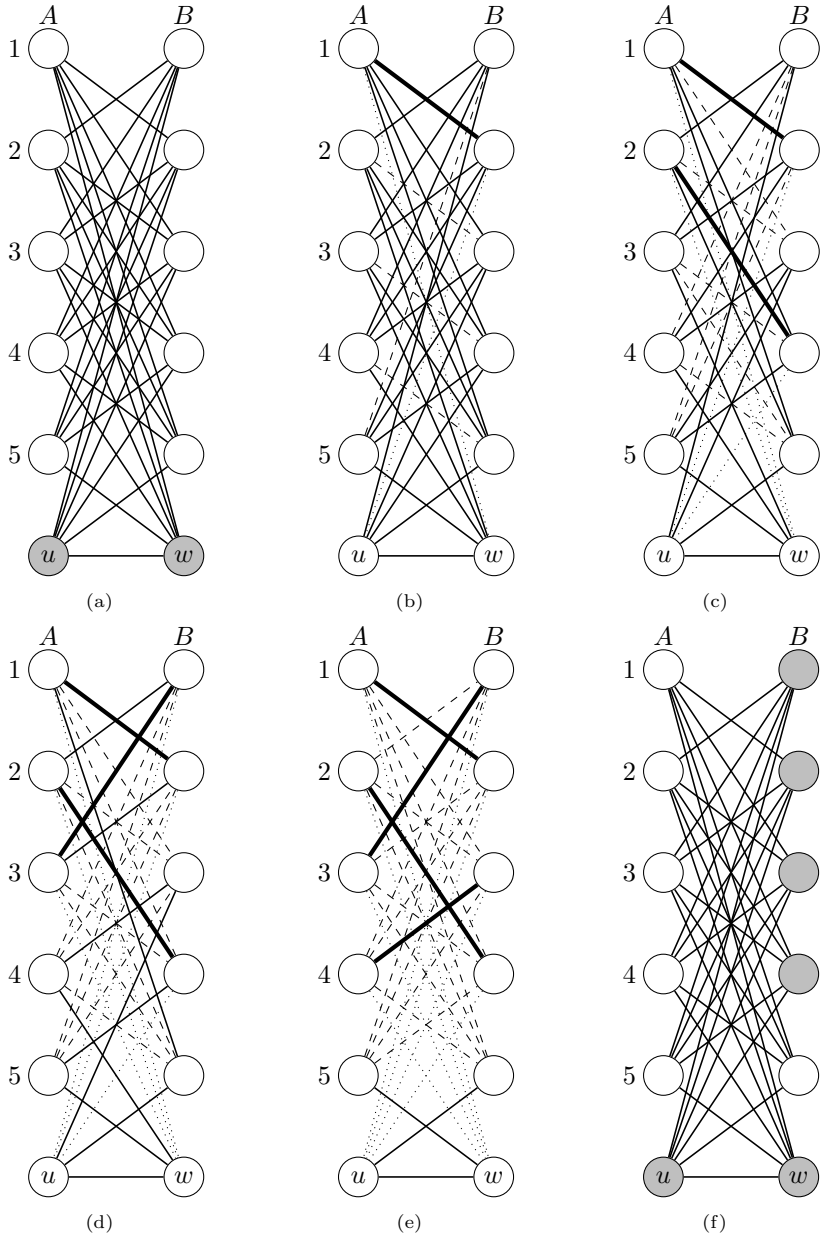

Figures~\ref{fig:tight_ex_3} and~\ref{fig:tight_ex_4} illustrate the process for $k=3$ and $k=4$, respectively. In these figures, the (perfect) matching identified by the algorithm at each iteration is depicted by dashed and dotted edges; the dashed lines identify the edges with both end-vertices in the optimal solution, and the dotted lines show the edges incident to $u$ or $w$. The edges in $\widetilde{E}$ are depicted by thick solid lines.
\begin{figure}[!ht]
\begin{subfigure}{.5\textwidth}
\centering
\begin{tikzpicture}[scale=0.45,auto,swap]
\foreach \pos/\name in {{(1,13)/A},{(5,13)/B},{(0,12)/1},{(0,9)/2},{(0,6)/3},{(0,3)/4}}
\node[label] (\name) at \pos {$\name$};
\foreach \pos/\name in {{(1,12)/A1}, {(5,12)/B1}, {(1,9)/A2}, {(5,9)/B2}, {(1,6)/A3}, {(5,6)/B3}, {(1,3)/A4}, {(5,3)/B4}}
\node[vertex] (\name) at \pos {};
\foreach \pos/\name in {{(1,0)/u}, {(5,0)/w}}
\node[darkVertex] (\name) at \pos {$\name$};
\foreach \source/ \dest /\weight in {A1/B2/, A1/B3/, A1/B4/, A2/B1/, A2/B3/, A2/B4/, A3/B1/, A3/B2/, A3/B4/, A4/B1/, A4/B2/, A4/B3/}
\path[semiEdge] (\source) -- node[weight]{}(\dest);
\foreach \source/ \dest /\weight in {A1/w/, A2/w/, A3/w/, A4/w/, B1/u/, B2/u/, B3/u/, B4/u/, u/w/}
\path[semiEdge] (\source) -- node[weight]{}(\dest);
\end{tikzpicture}
\caption{}
\label{fig:1ex_3}
\end{subfigure}
\begin{subfigure}{.5\textwidth}
\centering
\begin{tikzpicture}[scale=0.45,auto,swap]
\foreach \pos/\name in {{(1,13)/A},{(5,13)/B},{(0,12)/1},{(0,9)/2},{(0,6)/3},{(0,3)/4}}
\node[label] (\name) at \pos {$\name$};
\foreach \pos/\name in {{(1,12)/A1}, {(5,12)/B1}, {(1,9)/A2}, {(5,9)/B2}, {(1,6)/A3}, {(5,6)/B3}, {(1,3)/A4}, {(5,3)/B4}}
\node[vertex] (\name) at \pos {};
\foreach \pos/\name in {{(1,0)/u}, {(5,0)/w}}
\node[vertex] (\name) at \pos {$\name$};
\foreach \source/ \dest /\weight in {A1/B3/, A1/B4/, A2/B1/, A2/B4/, A3/B1/, A3/B2/, A4/B2/, A4/B3/}
\path[semiEdge] (\source) -- node[weight]{}(\dest);
\foreach \source/ \dest /\weight in {A2/w/, A3/w/, A4/w/, B1/u/, B3/u/, B4/u/, u/w/}
\path[semiEdge] (\source) -- node[weight]{}(\dest);
\foreach \source/ \dest /\weight in {A2/B3/, A3/B4/, A4/B1/}
\path[dashEdge] (\source) -- node[weight]{}(\dest);
\foreach \source/ \dest /\weight in {A1/B2/}
\path[ultraEdge] (\source) -- node[weight]{}(\dest);
\foreach \source/ \dest /\weight in {A1/w/, B2/u/}
\path[dotEdge] (\source) -- node[weight]{}(\dest);
\end{tikzpicture}
\caption{}
\label{fig:2ex_3}
\end{subfigure}
\begin{subfigure}{.32\textwidth}
\centering
\begin{tikzpicture}[scale=0.45,auto,swap]
\foreach \pos/\name in {{(1,13)/A},{(5,13)/B},{(0,12)/1},{(0,9)/2},{(0,6)/3},{(0,3)/4}}
\node[label] (\name) at \pos {$\name$};
\foreach \pos/\name in {{(1,12)/A1}, {(5,12)/B1}, {(1,9)/A2}, {(5,9)/B2}, {(1,6)/A3}, {(5,6)/B3}, {(1,3)/A4}, {(5,3)/B4}}
\node[vertex] (\name) at \pos {};
\foreach \pos/\name in {{(1,0)/u}, {(5,0)/w}}
\node[vertex] (\name) at \pos {$\name$};
\foreach \source/ \dest /\weight in {A1/B4/, A2/B1/, A3/B2/, A4/B3/}
\path[semiEdge] (\source) -- node[weight]{}(\dest);
\foreach \source/ \dest /\weight in {A3/w/, A4/w/, B1/u/, B3/u/, u/w/}
\path[semiEdge] (\source) -- node[weight]{}(\dest);
\foreach \source/ \dest /\weight in {A2/B3/, A3/B4/, A4/B1/, A1/B3/, A3/B1/, A4/B2/}
\path[dashEdge] (\source) -- node[weight]{}(\dest);
\foreach \source/ \dest /\weight in {A1/B2/, A2/B4/}
\path[ultraEdge] (\source) -- node[weight]{}(\dest);
\foreach \source/ \dest /\weight in {A1/w/, B2/u/, A2/w/, B4/u/}
\path[dotEdge] (\source) -- node[weight]{}(\dest);
\end{tikzpicture}
\caption{}
\label{fig:3ex_3}
\end{subfigure}
\begin{subfigure}{.33\textwidth}
\centering
\begin{tikzpicture}[scale=0.45,auto,swap]
\foreach \pos/\name in {{(1,13)/A},{(5,13)/B},{(0,12)/1},{(0,9)/2},{(0,6)/3},{(0,3)/4}}
\node[label] (\name) at \pos {$\name$};
\foreach \pos/\name in {{(1,12)/A1}, {(5,12)/B1}, {(1,9)/A2}, {(5,9)/B2}, {(1,6)/A3}, {(5,6)/B3}, {(1,3)/A4}, {(5,3)/B4}}
\node[vertex] (\name) at \pos {};
\foreach \pos/\name in {{(1,0)/u}, {(5,0)/w}}
\node[vertex] (\name) at \pos {$\name$};
\foreach \source/ \dest /\weight in {}
\path[semiEdge] (\source) -- node[weight]{}(\dest);
\foreach \source/ \dest /\weight in {A3/w/, B1/u/, u/w/}
\path[semiEdge] (\source) -- node[weight]{}(\dest);
\foreach \source/ \dest /\weight in {A2/B3/, A3/B4/, A4/B1/, A1/B3/, A3/B1/, A4/B2/, A1/B4/, A2/B1/, A3/B2/}
\path[dashEdge] (\source) -- node[weight]{}(\dest);
\foreach \source/ \dest /\weight in {A1/B2/, A2/B4/, A4/B3/}
\path[ultraEdge] (\source) -- node[weight]{}(\dest);
\foreach \source/ \dest /\weight in {A1/w/, B2/u/, A2/w/, B4/u/, A4/w/, B3/u/}
\path[dotEdge] (\source) -- node[weight]{}(\dest);
\end{tikzpicture}
\caption{}
\label{fig:4ex_3}
\end{subfigure}
\begin{subfigure}{.33\textwidth}
\centering
\begin{tikzpicture}[scale=0.45,auto,swap]
\foreach \pos/\name in {{(1,13)/A},{(5,13)/B},{(0,12)/1},{(0,9)/2},{(0,6)/3},{(0,3)/4}}
\node[label] (\name) at \pos {$\name$};
\foreach \pos/\name in {{(1,12)/A1}, {(5,12)/B1}, {(1,9)/A2}, {(1,6)/A3}, {(1,3)/A4}}
\node[vertex] (\name) at \pos {};
\foreach \pos/\name in {{(1,0)/u}, {(5,0)/w}}
\node[darkVertex] (\name) at \pos {$\name$};
\foreach \pos/\name in {{(5,9)/B2}, {(5,6)/B3}, {(5,3)/B4}}
\node[darkVertex] (\name) at \pos {};
\foreach \source/ \dest /\weight in {A1/B2/, A1/B3/, A1/B4/, A2/B1/, A2/B3/, A2/B4/, A3/B1/, A3/B2/, A3/B4/, A4/B1/, A4/B2/, A4/B3/}
\path[semiEdge] (\source) -- node[weight]{}(\dest);
\foreach \source/ \dest /\weight in {A1/w/, A2/w/, A3/w/, A4/w/, B1/u/, B2/u/, B3/u/, B4/u/, u/w/}
\path[semiEdge] (\source) -- node[weight]{}(\dest);
\end{tikzpicture}
\caption{}
\label{fig:5ex_3}
\end{subfigure}
\caption{Worst-case example for $k=3$: (a) the original graph; the set of white vertices is a maximum 3-dependent set, (b) the graph after removing $\mathcal{M}_1$, (c) the graph after removing $\mathcal{M}_2$, (d) the graph after removing $\mathcal{M}_3$, (e) the original graph; the gray vertices depict a minimum vertex cover of the residual graph, i.e., after removing $\mathcal{M}_3$, and the set of white vertices is the output of the algorithm.}
\label{fig:tight_ex_3}
\end{figure}
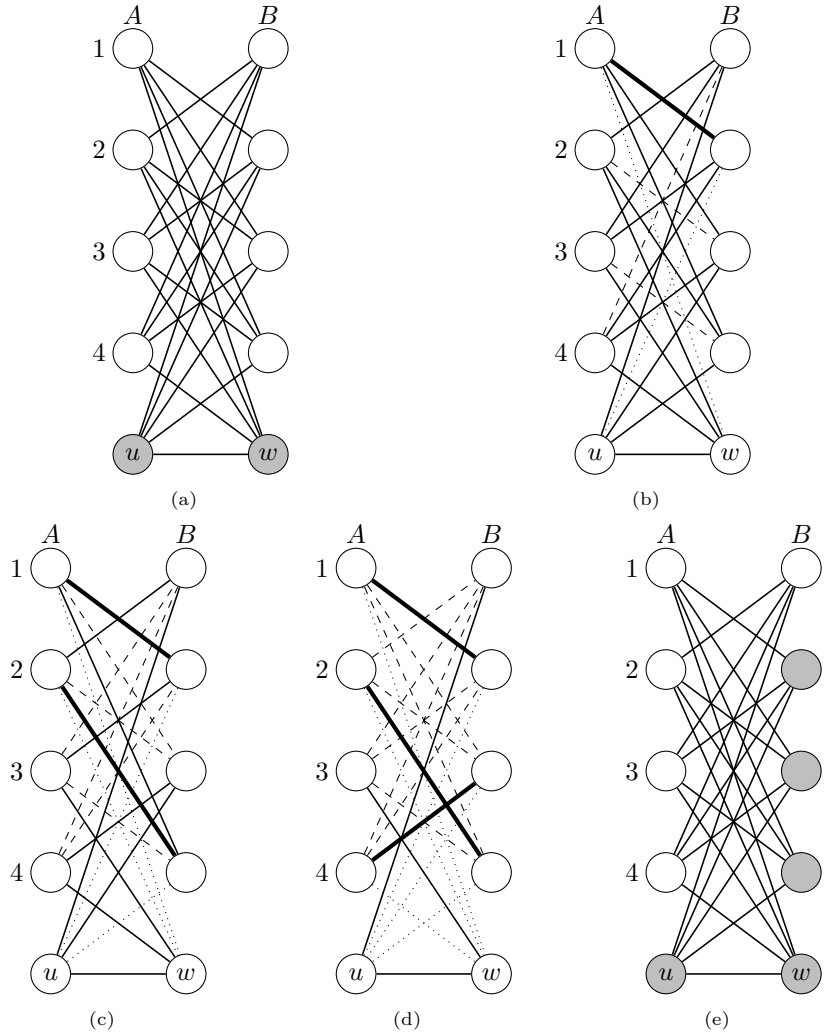
%
\section{Conclusion} \label{sec:conclusion}
\noindent
The {\sc Maximum $k$-dependent Set} problem is NP-hard on the class of bipartite graphs, for any $k \geq 1$. This paper presents a $(1+\frac{k}{k+2})$-approximation algorithm for this problem that runs in $\mathcal{O}(km \sqrt{n})$ time on a (bipartite) graph with $n$ vertices and $m$ edges and improves upon its previously best-known approximation ratio of $1+\frac{k}{k+1}$. Our proof also indicates that the algorithm retains its approximation ratio when applied to K\"{o}nig-Egerv\'{a}ry (KE) graphs, which is a supperclass of bipartite graphs. Furthermore, we present a family of bipartite graphs, on which the algorithm achieves the approximation bound.

As an extension of this work, one may replace the loop of Algorithm~\ref{algorithm}, i.e., removing $k$ maximum-cardinality matchings, with removing a single maximum-cardinality $k$-matching and study the approximation ratio of the new algorithm. A $k$-matching in a graph is a subset of edges such that every vertex in the graph is incident to at most $k$ of them. The maximum $k$-matching problem, which is to find a $k$-matching with the largest number of edges in a given graph, is polynomial-time solvable~\cite{Gerards95}. The structure of the examples presented in this paper suggests such an algorithm may avoid the pitfall of the current method in missing the ``correct'' edges, hence achieving an approximation ratio better than $1+\frac{k}{k+2}$. In the presence of multiple maximum-cardinality matchings in a graph, one may also try to quantify the probability of missing the ``correct'' edges---based on the employed matching subroutine---which may lead to the design of a randomized variant of the current algorithm with a better expected approximation ratio. The proposed algorithm always generates a feasible solution for the {\sc Maximum $k$-dependent Set} problem on general graphs. Also, the minimum cardinality of a vertex cover is an upper bound on the maximum cardinality of a matching in every graph. Let $\theta(G)$ denote the gap between these two invariants in a graph $G$. Then, for the classes of graphs with polynomial-time computable independence number in the residual, it is interesting to establish the approximation ratio of the algorithm presented in this paper as a function of $\theta(G)$ and $k$.

\paragraph*{Acknowledgements}
Partial support by AFOSR award FA9550-19-1-0161 is  gratefully acknowledged. {We would also like to thank the four anonymous referees for their insightful comments and suggestions.}


\end{document}